\documentclass[12pt,a4paper]{article}

\usepackage{amsmath,amssymb,amsfonts,amsthm,graphicx}

\newtheorem{theorem}{Theorem}

\newtheorem{prop}[theorem]{Proposition}
\newtheorem{problem}[theorem]{Problem}

\newcommand{\Ga}{\ensuremath{\Gamma}}

\newcommand{\diam}{\rm diam}

\begin{document}

\title{Eigenvalues of graphs and spectral Moore theorems}
\author{‹ž"Sebastian M. Cioab\u{a}\footnote{Department of  Mathematical Sciences, University of Delaware, Newark, DE 19716-2553, USA. This research has been partially supported by NSF grants DMS-1600768 and CIF-1815922 and a and a JSPS Invitational Fellowship for Research in Japan S19016.}}
\date{\today}
\maketitle

\begin{abstract}
In this paper, we describe some recent spectral Moore theorems related to determining the maximum order of a connected graph of given valency and second eigenvalue. We show how these spectral Moore theorems have applications in Alon-Boppana theorems for regular graphs and in the classical degree-diameter/Moore problem. 
\end{abstract}

\section*{\S1 Introduction}

Our graph theoretic notation is standard (see \cite{BH,Bollobas}). Let $\Ga=(V,E)$ be an undirected graph with vertex set $V$ and edge set $E$. Given $u,v\in V$, the distance $d(u,v)$ equals the minimum length of a path between $u$ and $v$ if such a path exists or $\infty$ otherwise. If $\Ga$ is connected, then all the distances between its vertices are finite and the diameter $\diam(\Ga)$ of $\Ga$ is defined as the maximum of $d(u,v)$, where the maximum is taken over all pairs $u,v\in V$. The pairwise distance and the diameter of a connected graph can be calculated efficiently using breadth first search. The Moore or degree-diameter problem is a classical problem in combinatorics (see \cite{MS}).
\begin{problem}\label{mooreprob}
Given $r\geq 3$ and $D\geq 2$, what is the maximum order $n_{r,D}$ of a connected $r$-regular graph of diameter $D$ ?
\end{problem}
There is a well known upper bound for $n_{r,D}$ known as the Moore bound $m_{r,D}$ which is obtained as follows. If $\Ga$ is a connected $r$-regular graph of diameter $D$, then for any given vertex $x$ in $\Ga$ and any $1\leq j\leq D$, the number of vertices at distance $j$ from $x$ is at most $r(r-1)^{j-1}$. Therefore, 
\begin{equation}\label{moorebnd}
n_{r,D}\leq 1+r+r(r-1)+\dots +r(r-1)^{D-1}.
\end{equation}
We denote by $m_{r,D}$ the right hand-side of the above inequality. When $r=2$, it is straightforward to note that $n_{2,D}=m_{2,D}=2D+1$ and the maximum is attained by the cycle $C_{2D+1}$ on $2D+1$ vertices.  For $D=1$, it is easy to see that $n_{r,1}=m_{r,1}=r+1$ and the maximum is attained by the complete graph $K_{r+1}$ on $r+1$ vertices. 

For $D=2$, a classical result of Hoffman and Singleton \cite{HS} gives that $n_{r,2}$ equals the Moore bound $m_{r,2}=r^2+1$ only when $r=2$ (attained by the cycle $C_5$), $r=3$ (the Petersen graph), $r=7$ (the Hoffman-Singleton graph) or possibly $r=57$. The existence of a $57$-regular graph with diameter $2$ on $57^2+1=3250$ vertices is a well known open problem in this area (see \cite{Dalfo,JV,MaSi}). For $D\geq 3$ and $r\geq 3$, Damerell \cite{D} and independently, Bannai and Ito \cite{BI} proved that there are no graphs attaining the Moore bound \eqref{moorebnd}.

The adjacency matrix $A$ is the $V\times V$ matrix whose $(x,y)$-th entry equals the number of edges between $x$ and $y$. This matrix is a real symmetric matrix and if $\Ga$ is simple (no loops nor multiple edges), then $A$ is a $(0,1)$ symmetric matrix.  Let $r\geq 3$ be a given integer. We will use the following family of orthogonal polynomials:
\begin{align}
F_0(x)&=1, F_1(x)=x, F_2(x)=x^2-r, \label{polyf1}\\
F_{j}(x)&=xF_{j-1}(x)-(r-1)F_{j-2}(x), \label{polyf2}
\end{align}
for any $j\geq 3$. Let $q=\sqrt{k-1}$. The polynomials $(F_i)_{i\geq 0}$ form a sequence of orthogonal polynomials 
with respect to the positive weight
\[
w(x)=\frac{\sqrt{4q^2-x^2}}{k^2-x^2}
\]
on the interval $[-2q, 2q]$  (see \cite[Section 4]{HOb}). The polynomials $F_i(q y)/ q^i$ in $y$ are called Geronimus polynomials \cite{Ger}. 

For any vertices $u$ and $v$ of $\Ga$ and any non-negative integer $\ell$, the entry $(u,v)$ of the matrix $A^{\ell}$ equals the number of walks of length $\ell$ between $u$ and $v$. A walk $u=u_0,u_1,\dots,u_{\ell-1},u_{\ell}=v$ in $G$ is called non-backtracking if $u_{i}u_{i+1}\in E$ for any $0\leq i\leq \ell-1$ and $u_{i}\neq u_{i+2}$ for any $0\leq i\leq \ell-2$ (when $\ell\geq 2$). The  following result goes back to Singleton \cite{Singleton}.
\begin{prop}[Singleton \cite{Singleton}]
Let $\Ga$ be a connected $r$-regular graph with adjacency matrix $A$. For any vertices $u$ and $v$ of $\Ga$ and any non-negative integer $\ell$, the entry $(u,v)$ of the matrix $F_{\ell}(A)$ equals the number of non-backtracking walks of length $\ell$ between $u$ and $v$.
\end{prop}
The eigenvalues of $A$ are real and we denote them by $\lambda_1\geq \lambda_2\geq \dots \geq \lambda_n$, where $n=|V|$. Sometimes, to highlight the dependence of the eigenvalues on a particular graph, we will use $\lambda_j(\Ga)$ for $\lambda_j$. When $\Ga$ is $r$-regular and connected, it is known that $\lambda_1=r$ and that $\lambda_2<r$. It is also known that $\lambda_j\in [-r,r]$ and that $-r$ is an eigenvalue if and only if the graph is bipartite. The smallest eigenvalue of a regular graph has been used to determine the independence number of various interesting graphs (see Godsil and Meagher \cite{GM}).

The properties of the eigenvalues of a regular graph were essential in the proofs of Hoffman and Singleton \cite{HS} as well as Damerell \cite{D} and Bannai and Ito \cite{BI}. 

The spectral gap $r-\lambda_2$ is an important parameter in spectral graph theory and is closely related to the connectivity \cite{Fiedler} and expansion properties of the graph \cite{HLW}. Informally, expanders are sparse graphs with large spectral gap. More precisely, a family $(\Ga_{m})_{m\geq 1}$ of graphs is called a family of expanders if 
\begin{enumerate}
\item there exists $r\geq 3$ such that each $\Ga_m$ is a connected $r$-regular graph for $m\geq 1$ and the number of vertices of $\Ga_m$ goes to infinity as $m$ goes to infinity,
\item there is a positive constant $c_r>0$ such that $r-\lambda_2(\Ga_m)>c_r$ for any $\Ga_m$.
\end{enumerate}
The first condition above explains the denomination of {\em sparse} used at the beginning of this paragraph. This condition implies that the number of edges in $\Ga_m$ is linear in its number of vertices for every $m\geq 1$. This is best possible in order of magnitude for connected graphs. 

The second condition is algebraic and is equivalent to a combinatorial condition that each $\Ga_m$ is highly connected (meaning their expansion constants are bounded away from $0$) and also equivalent to the probability condition that a random walk on $\Ga_m$ converges quickly to its stationary distribution. We refer to \cite{HLW} for the precise descriptions of these conditions.

A natural question arising from the previous considerations is how large can the spectral gap $r-\lambda_2(\Ga)$ be for an $r$-regular connected graph $\Ga$ ? Since we are interested in situations where $r\geq 3$ is fixed, this is equivalent to asking how small can $\lambda_2$ be for a connected $r$-regular graph $\Ga$. This was answered by the Alon-Boppana theorem.

\begin{theorem}
Let $r\geq 3$ be a natural number.
\begin{enumerate}
\item {\bf Alon-Boppana 1986.} If $\Ga$ is a connected $r$-regular graph with $n$ vertices, then
\begin{equation}
\lambda_2(\Ga)\geq 2\sqrt{r-1}\left(1-\frac{C}{\diam^2(\Ga)}\right)=2\sqrt{r-1}(1-o(1)),
\end{equation}
where $C>0$ is a constant and $o(1)$ is a quantity that goes to $0$ as $n$ goes to infinity.
\item {\bf Asymptotic Alon-Boppana Theorem.} If $(\Ga_m)_{m\geq 1}$ is a sequence of connected $r$-regular graphs such that $|V(\Ga_m)|\rightarrow \infty$ as $m\rightarrow \infty$. Then
\begin{equation}
\liminf_{m\rightarrow \infty}\lambda_2(\Ga_m)\geq 2\sqrt{r-1}.
\end{equation}
\end{enumerate}
\end{theorem}

To my knowledge, there is no paper written by Alon and Boppana which contains the theorem above. The first appearance of this result that I am aware of, is in 1986 in Alon's paper \cite{A1} where it is stated that
\begin{center}
{\em R. Boppana and the present author showed that for every $d$-regular graph $G$ on $n$ vertices $\lambda(G)\leq d-2\sqrt{d-1}+O(\log_d n)^{-1}$.} 
\end{center}
Note that the $\lambda$ in \cite{A1} is the smallest positive eigenvalue of the Laplacian $D-A=rI-A$ of $G$ and it equals $r-\lambda_2(G)$. Therefore the statement above is equivalent to 
\begin{equation}
\lambda_2(G)\geq 2\sqrt{r-1}-O(\log_r n)^{-1}.
\end{equation}
There is a similar result to the Alon-Boppana theorem that is due to Serre \cite{Serre}. The meaning of this theorem below is that large $r$-regular graphs tend to have a positive proportion of eigenvalues trying to be greater than $2\sqrt{r-1}$.
\begin{theorem}[Serre \cite{Serre}]\label{serrethm}
For any $r\geq 3,\epsilon>0$, there exists $c=c(\epsilon,r)>0$ such that any $r$-regular graph $\Ga$ on $n$ vertices has at least $c\cdot n$ eigenvalues that are at least $2\sqrt{r-1}-\epsilon$.
\end{theorem}
These results motivated the definition of Ramanujan graphs that was introduced by Lubotzky, Phillips and Sarnak \cite{LPS}. A connected $r$-regular graph $\Ga$ is called Ramanujan if all its eigenvalues (with the exception of $r$ and perhaps $-r$, if $\Ga$ is bipartite) have absolute value at most $2\sqrt{r-1}$. Lubotzky, Phillips and Sarnak \cite{LPS} and independently Margulis \cite{Margulis} constructed infinite families of $r$-regular Ramanujan graphs when $r-1$ is a prime. These constructions used results from algebra and number theory closely related to a conjecture of Ramanujan regarding the number of ways of writing a natural number as a sum of four squares of a certain kind (see \cite{LPS,Margulis} and also, \cite{DSV} for a more detailed description of these results). For the longest time, it was not known whether infinite families of $r$-regular Ramanujan graphs exist for any $r\geq 3$. Marcus, Spielman and Srivastava \cite{MSS} obtained a breakthrough result by showing that there exist infinite families of bipartite $r$-regular Ramanujan graphs for any $r\geq 3$. Their method of interlacing polynomials has been fundamental to this proof and has found applications in other areas of mathematics as well.

\section*{\S2 Spectral Moore theorems for general graphs}

Throughout the years, several proofs of the Alon-Boppana theorem have appeared (see Lubotzky, Phillips and Sarnak \cite{LPS}, Nilli \cite{Nilli1}, Kahale \cite{Kahale}, Friedman \cite{Fri1}, Feng and Li \cite{FengLi}, Li and Sol\'{e} \cite{LiSole}, Nilli \cite{Nilli2} and Mohar \cite{Moh}). Theorem \ref{serrethm} was proved by Serre \cite{Serre} with a non-elementary proof (see also \cite{DSV}). The first elementary proofs appeared around the same time by Cioab\u{a} \cite{Ci06} and Nilli \cite{Nilli2}. Richey, Stover and Shutty \cite{RSS} worked to turn Serre's proof into a quantitative theorem and asked the following natural question.
\begin{problem}\label{vrt}
Given an integer $r\geq 3$ and $\theta<2\sqrt{r-1}$, what is the maximum order $v(r,\theta)$ of a $r$-regular graph $\Ga$ with $\lambda_2(\Ga)\leq \theta$ ?
\end{problem}
These authors obtained several results involving $v(r,\theta)$. In this section, we describe our recent results related to the problem above and its bipartite and hypergraph versions. See \cite{CKN19,CKNV16,CKMNO20} and the references therein for more details and other related problems. The method that is fundamental to all these results is due to Nozaki \cite{NozakiLP} who proved the {\em linear programming bound for graphs}.
\begin{theorem}[Nozaki \cite{NozakiLP}]\label{nozakithm}
Let $\Ga$ be a connected $r$-regular graph with $v$ vertices and distinct eigenvalues $\theta_1=k>\theta_2> \ldots >\theta_d$.  If there exists a polynomial $f(x)=\sum_{i=0}^{t}  f_i F_i(x)$ such that  
$f(r)>0$, $f(\theta_i) \leq 0$ for any $2\leq i\leq d$, $f_0>0$, and $f_i \geq 0$ for any $1\leq i\leq t$, then
\begin{equation*} 
v \leq \frac{f(r)}{f_0}. 
\end{equation*}
\end{theorem}
Nozaki used this result to study the following problem.
\begin{problem}\label{nozakiproblem}
Given integers $v>r\geq 3$, what is the $r$-regular graph $\Ga$ on $v$ vertices that has the smallest $\lambda_2$ among all $r$-regular graphs on $v$ vertices ?
\end{problem}
While similar to it, this problem is quite different from Problem \ref{vrt}. 

In \cite{CKNV16}, the authors used Nozaki's LP bound for graphs to obtain the following general upper bound for $v(r,\theta)$.
\begin{theorem}[Cioab\u{a}, Koolen, Nozaki and Vermette \cite{CKNV16}]\label{cknv16thm}
Given integers $r,t\geq 3$ and a non-negative real number $c$, let $T(r,t,c)$ be the $t \times t$ tridiagonal matrix: 
{\small
$$
T(r,t,c)= \left[\begin{array}{ccccccc}
0 & r\\
1 & 0 & r-1\\
    & 1 & 0 & r-1 \\
    &     &  .  &  . & . \\
    &     &     &  . & . & . \\
    &     &     &    & 1 & 0 & r-1 \\
    &     &     &    &         &  c   & r-c
\end{array} \right]
$$
}
If $\theta$ equals the second largest eigenvalue $\lambda_2(T(r,t,c))$ of the matrix $T(r,t,c)$, then
$$
v(r,\theta) \leq 1+\sum_{i=0}^{t-3} r(r-1)^i+ \frac{r(r-1)^{t-2}}{c}. 
$$
\end{theorem}

We sketch below the ideas of the proof of this theorem. For $j\geq 0$, denote 
\begin{equation}\label{polyg}
G_{j}=\sum_{i=0}^{j}F_i, 
\end{equation}
where the $F_i$s are the orthogonal polynomials defined in equations \eqref{polyf1} and \eqref{polyf2}. The polynomials $(G_j)_{j\geq 0}$ also form a family of orthogonal polynomials. They satisfy the following properties:
\begin{align}
G_0(x)&=1, G_1(x)=x+1, G_2(x)=x^2+x-(r-1) \label{polyg1}\\
G_j(x)&=xG_{j-1}-(r-1)G_{j-2}(x), \label{polyg2}
\end{align}
for $j\geq 3$.

The eigenvalues of the matrix $T=T(r,t,c)$ are the roots of $(x-r)(G_{t-1}+(c-1)G_{t-2})$ and are distinct (see \cite[Theorem 2.3]{CKNV16}). If we denote them by $r=\lambda_1>\lambda_2>\dots >\lambda_t$, then the polynomial $f(x)=\frac{1}{c}\cdot (x-\lambda_2)\prod_{i\geq 3}(x-\lambda_i)^2$ satisfies $f(\lambda_i)\leq 0$ for $i\geq 2$. It is a bit more involved to check the other conditions from Theorem \ref{nozakithm} and we refer the reader to \cite{CKNV16} for the details to see how one can apply Nozaki's LP bound to $f$ and obtain that
$$
v(r,\theta)\leq \frac{f(r)}{f_0}=\sum_{i=0}^{t-2}F_i(r)+F_{t-1}(r)/c=1+\sum_{i=0}^{t-3} r(r-1)^i+ \frac{r(r-1)^{t-2}}{c}.
$$

To make things more clear, note the following result.
\begin{prop}
Let $r\geq 3$ be an integer. For any $\theta\in [-1,2\sqrt{r-1})$, there exists an integer $t$ and a positive number $c$ such that $\theta$ is the second largest eigenvalue of the matrix $M(r,t,c)$. 
\end{prop}
Let $\lambda^{(t)}$ denote the largest root of $G_t$ and $\mu^{(t)}$ denote the largest root of $F_t$. Note that $\lambda^{(1)}=-1<0=\mu^{(1)}$ and 
\begin{equation}\label{lambda2mu2}
\lambda^{(2)}=\frac{-1+\sqrt{4r-3}}{2}<\sqrt{r}=\mu^{(2)}.
\end{equation}
Bannai and Ito \cite[Section III.3]{BIBook} showed that $\lambda^{(t)}=2\sqrt{r-1}\cos\tau$, where $\frac{\pi}{t+1}<\tau<\frac{\pi}{t}$. Because $F_t=G_t-G_{t-1}$, one can show that $\lambda^{(t)}<\mu^{(t)}$ for any $t$ (see \cite[Prop 2.6]{CKNV16} for other properties of these eigenvalues). From the remarks following Theorem \ref{cknv16thm}, note that the second largest eigenvalue $\lambda_2(t,c)$ of $T(r,t,c)$ equals the largest root of the polynomial $(c-1)G_{t-1}+G_{t-2}$. Because the roots of $G_{t-2}$ and $G_{t-1}$ interlace, one obtains that $\lambda_2(t,c)$ is a decreasing function in $c$ and takes values between $\lim_{c\rightarrow \infty}\lambda_2(t,c)=\lambda^{(t-2)}$ and $\lim_{c\rightarrow 0}\lambda_2(t,c)=\mu^{(t-1)}$. Taking into the account what happens when $c=1$, namely that $\lambda_2(t,c)=\lambda^{(t-1)}$, we obtain the following result.
\begin{prop}\label{lambda2tc}
For $c\in [1,\infty)$, $\lambda_2(t,c)$ takes any value in the interval $[\lambda^{(t-1)},\lambda^{(t-2)})$. 
\end{prop}
Putting these things together, one deduces that $\lambda_{2}(t,c)$ can take any possible value between $\lambda_2(2,1)=-1$ and $\lim_{t\rightarrow \infty}\lambda_2(t,c)=2\sqrt{r-1}$. There are several infinite families $(r,\theta)$ for which the precise values $v(r,\theta)$ have been determined in \cite{CKNV16}, but there are several open problems for relatively small values of $r$ and $\theta$. For example, $v(6,2)\geq 42$ with an example of a $6$-regular graph with $\lambda_2=2$ on $42$ vertices being the 2nd subconstituent of the Hoffman-Singleton graph. Theorem \ref{cknv16thm} can give $v(6,2)\leq 45$ (see \cite{CKMNO20}). Also, we know that $v(3,\sqrt{2})=14$ (Heawood graph), but we don't know the exact value of $v(k,\sqrt{2})$ for any $k\geq 3$. Lastly, $v(k,\sqrt{k})$ has been determined for $k=3$ (equals $18$ with Pappus graph as an example attaining it) and $k=4$ (it is $35$ with the odd graph $O_4$ meeting it), but we don't know it for $k\geq 5$.

\section*{\S3 Alon-Boppana and Serre theorems}

We point out the relevance of these results in the context of Alon-Boppana and Serre theorems. A typical Alon-Boppana result is of the form: if $\Ga$ is a connected $r$-regular graph with diameter $D\geq 2k$, then 
\begin{equation}\label{ab}
\lambda_2(\Ga)\geq 2\sqrt{r-1}\cos\frac{\pi}{k+1}.
\end{equation}
See Friedman \cite[Corollary 3.6]{Fri1} for the inequality above or Nilli \cite[Theorem 1]{Nilli2} for a slightly weaker bound. The equivalent contrapositive formulation of inequality \eqref{ab} is the following: if $\Ga$ is a connected $r$-regular graph of diameter $D$ with $\lambda_2(\Ga)<2\sqrt{r-1}\cos\frac{\pi}{k+1}$, then $D<2k$. By Moore bound \eqref{moorebnd}, this implies that 
\begin{equation}\label{ab1}
|V(\Ga)|\leq m_{r,2k-1}=1+r+r(r-1)+\dots +r(r-1)^{2k-2}.
\end{equation}

Obviously, the best bound one can achieve here is obtained for that $k$ with $2\sqrt{r-1}\cos\frac{\pi}{k}\leq \lambda_2(\Ga)<2\sqrt{r-1}\cos\frac{\pi}{k+1}$. When applying Theorem \ref{cknv16thm}, let $\theta$ be a real number such that $2\sqrt{r-1}\cos\frac{\pi}{k}\leq \theta <2\sqrt{r-1}\cos\frac{\pi}{k+1}$. Given the properties of the largest roots $\lambda^{(t)}$ of the polynomials $G_t$ (see Proposition \ref{lambda2tc} and the paragraph containing it), we must have that either $\theta\in (\lambda^{(k-1)},\lambda^{(k)}]$ or $\theta \in (\lambda^{(k)},\lambda^{(k+1)})$. If $\theta\in (\lambda^{(k-1)},\lambda^{(k)}]$, then there exists $c_1\geq 1$ such that $\theta$ is the largest eigenvalue of the polynomial $G_{k}+(c_1-1)G_{k-1}$. If $\Ga$ is a connected $r$-regular graph with $\lambda_2(\Ga)\leq \theta$, then using Theorem \ref{cknv16thm} we get that
\begin{equation}\label{ourab1}
|V(\Ga)|\leq v(r,\theta)=1+r+r(r-1)+\dots +r(r-1)^{k-2}+\frac{r(r-1)^{k-1}}{c_1}
\end{equation}
which is clearly better than \eqref{ab1}. If $\theta\in (\lambda^{(k)},\lambda^{(k+1)})$, then there exists $c_2>1$ such that $\theta$ is the largest eigenvalue of the polynomial $G_{k+1}+(c_2-1)G_{k}$. As above, if $\Ga$ is a connected $r$-regular graph with $\lambda_2(\Ga)\leq \theta$, then Theorem \ref{cknv16thm} implies that
\begin{equation}\label{ourab2}
|V(\Ga)|\leq v(r,\theta)=1+r+r(r-1)+\dots +r(r-1)^{k-1}+\frac{r(r-1)^k}{c_2}
\end{equation}
which is again better than \eqref{ab}. 


\section*{\S4 Spectral Moore theorems for bipartite graphs}

Building on this work, the author with Koolen and Nozaki extended and refined these results to bipartite regular graphs \cite{CKN19}. Let $r\geq 3$ be an integer and $\theta$ be any real number between $0$ and $2\sqrt{r-1}$. Define $b(r,\theta)$ as the maximum number of vertices of a bipartite $r$-regular graph whose second largest eigenvalue is at most $\theta$. Clearly, $b(r,\theta)\leq v(r,\theta)$ and a natural question is whether or not these parameters are actually the same or not. One can show that $v(3,1)=10$ attained by the Petersen graph and that $b(3,1)=8$ attained by the $3$-dimensional cube. For the bipartite graphs, a linear programming bound similar to Nozaki's Theorem \ref{nozakithm} from \cite{NozakiLP} was obtained with the use of the following polynomials:
\begin{equation*}
\mathcal{F}_{0,i}=F_{2i}(\sqrt{x}), \mathcal{F}_{1,i}=\frac{F_{2i+1}(\sqrt{x})}{\sqrt{x}},
\end{equation*}
for any $i\geq 0$, where $(F_i)_{i\geq 0}$ were defined earlier in \eqref{polyf1} and \eqref{polyf2}. Let $\Ga$ be a bipartite connected regular graph. Its adjacency matrix $A$ has the form $\begin{bmatrix} 0 & \boldsymbol{N}\\\boldsymbol{N}^{\top} & 0\end{bmatrix}$ and we call $\boldsymbol{N}$ the biadjacency matrix of $\Ga$. Note that $F_{2i}(A)=\begin{bmatrix}\mathcal{F}_{0,i}(\boldsymbol{N}\boldsymbol{N}^{\top}) & 0\\0 & \mathcal{F}_{0,i}(\boldsymbol{N}^{\top}\boldsymbol{N})\end{bmatrix}$ for any $i\geq 0$. The following is called the LP bound for bipartite regular graphs.
\begin{theorem}[Cioab\u{a}, Koolen and Nozaki \cite{CKN19}]\label{cknLPthm}
Let $\Gamma$ be a connected bipartite $r$-regular graph with $v$ vertices and denote by $\{\pm \tau_0, \ldots, \pm \tau_d \}$ its set of distinct eigenvalues, where $\tau_0=r$. If there exists a polynomial $f(x)=\sum_{i= 0}^t  f_i \mathcal{F}_{0,i}(x)$ such that 
$f(r^2)>0$, $f(\tau_i^2) \leq 0$ for each $i\in  \{1,\ldots, d \}$, $f_0>0$, and $f_j \geq 0$ for each $j \in \{1,\ldots, t \}$, then 
\begin{equation} 
v \leq  \frac{2 f(r^2)}{f_0}. 
\end{equation}
Equality holds if and only if for each $i \in \{1,\ldots, d\}$, $f(\tau_i^2)=0$ and for each $j\in \{1,\dots,t\}$, ${\rm tr}(f_j\mathcal{F}_{0,j}(\boldsymbol{N}\boldsymbol{N}^{\top}))=0$, 
and ${\rm tr} (f_j\mathcal{F}_{0,j}(\boldsymbol{N}^{\top} \boldsymbol{N}))=0$, 
where $\boldsymbol{N}$ is the biadjacency matrix of $\Gamma$.
If equality holds and $f_j>0$ for each $j \in \{1,\ldots, t\}$, then the girth of $\Gamma$ is at least  
$2t+2$. 
\end{theorem}
For any integers $t\geq 3, r\geq 3$ and any positive $c\leq r$, let $\boldsymbol{B}(r,t,c)$ be the $t\times t$ tridiagonal matrix with lower diagonal $(1,\ldots,1,c,r)$, upper diagonal $(r,r-1,\ldots, r-1,r-c)$, and constant row sum $r$. Using Theorem \ref{cknLPthm}, the following general upper bound for $b(r,\theta)$ was obtained in \cite{CKN19}.
\begin{theorem}[Cioab\u{a}, Koolen and Nozaki \cite{CKN19}]\label{ckn19thm}
If $\theta$ is the second largest eigenvalue of $\boldsymbol{B}(r,t,c)$, then
\begin{equation}
b(r,\theta) \leq 2\left(\sum_{i=0}^{t-4} (r-1)^{i}+ \frac{(r-1)^{t-3}}{c}+\frac{(r-1)^{t-2}}{c}\right):=M(r,t,c). 
\end{equation}
Equality holds if and only if there exists a bipartite distance-regular graph whose quotient matrix with respect to the distance-partition from a vertex is $\boldsymbol{B}(r,t,c)$ for $1\leq c<r$ or $\boldsymbol{B}(r,t-1,1)$ for $c=r$.   
\end{theorem} 
Define $H_j(x)=\sum_{i=0}^{\lfloor j/2\rfloor}F_{j-2i}(x)$ for $j\geq 0$. These are orthogonal polynomials and one can show that $H_j(x)=xH_{j-1}(x)-(r-1)H_{j-2}(x)$ for $j\geq 2$ as well as that $H_{j}(x)=\frac{F_{j+2}(x)-(r-1)^2F_{i}(x)}{x^2-r^2}$. The first step in proving the above result is showing that the characteristic polynomial of $\boldsymbol{B}(r,t,c)$ equals $(x^2-r^2)(H_{t-2}(x)+(c-1)H_{t-4}(x))$. The proof proceeds in similar steps to Theorem \ref{cknv16thm}, but is more technical and we refer the reader to \cite{CKN19} for the details. Similar to the situation for general graphs, one can show the following.
\begin{prop}
Let $r\geq 3$ be an integer. For any $\theta\in [0,2\sqrt{r-1})$, there exists $t$ and $c$ such that $\theta$ is the second largest eigenvalue of $\boldsymbol{B}(r,t,c)$.
\end{prop}
In \cite{CKN19}, the authors also proved that for given $r$ and $\theta$, the upper bound obtained in Theorem \ref{ckn19thm} is better than the one in Theorem \ref{cknv16thm}. Theorem \ref{ckn19thm} has applications to various areas and it improves results obtained in the context of coding theory by H\o holdt and Janwa \cite{HJa} and H\o holdt and Justensen \cite{HJ} and design theory by Teranishi and Yasuno \cite{TY}. 

As in the case of Theorem \ref{cknv16thm}, Theorem \ref{ckn19thm} has applications for the Alon-Boppana theorems for bipartite regular graphs. Corollary 4.11 in \cite{CKN19} is a consequence of Theorem \ref{ckn19thm} and states that if $\Ga$ is a bipartite $r$-regular graph of order greater than $M(r,t,c)$ (the right hand-side in Theorem \ref{ckn19thm}), then $\lambda_2\geq \theta$, where $\theta$ is the second largest eigenvalue of $\boldsymbol{B}(r,t,c)$. Li and Sol\'{e} \cite[Theorems 3 and 5]{LiSole} proved that if $\Ga$ is a bipartite $r$-regular graph of girth $2\ell$, then $\lambda_2(\Ga)\geq 2\sqrt{r-1}\cos\frac{\pi}{\ell}$. This result follows from Corollary 4.11 in \cite{CKN19} as $2\sqrt{r-1}\cos\frac{\pi}{\ell}$ is the second largest eigenvalue of $\boldsymbol{B}(r,\ell+1,1)$ and having girth $2\ell$ implies that $\Ga$ has at least $M(r,\ell+1,1)$ vertices. 

\section*{\S5 Classical Moore problem}

Since the fundamental work of Singleton \cite{Singleton}, Hoffman and Singleton \cite{HS}, Bannai and Ito \cite{BI} and Damerell \cite{D} in the 1970s, the families of orthogonal polynomials $(F_j)_{j\geq 0}$ and $(G_{j})_{j\geq 0}$ have been important in the study of the Moore problem \eqref{mooreprob}. It has been observed by several authors (see \cite{DSS} or \cite{MS} for example) that if $\Ga$ is connected $r$-regular of diameter $D$, eigenvalues $r=\lambda_1>\lambda_2\geq \dots \geq \lambda_n$ and $\beta=\max(|\lambda_2|,|\lambda_n|)$, then 
\begin{equation}\label{moorebndlarge}
|V(\Ga)|\leq G_D(r)-G_D(\beta)=m_{r,D}-G_D(\beta),
\end{equation}
where $m_{r,D}$ is the upper bound from the Moore bound \eqref{moorebnd}. Recall that $\lambda^{(D)}$ denotes the largest root of $G_D(x)$ and satisfies 
\begin{equation}
2\sqrt{r-1}\cos\frac{\pi}{D}<\lambda^{(D)}<2\sqrt{r-1}\cos\frac{\pi}{D+1}.
\end{equation}
Inequality \eqref{moorebndlarge} will improve the classical Moore bound \eqref{moorebnd} when $G_D(\beta)>0$. This will happen when $\beta>\lambda^{(D)}$. When $D=2$, from \eqref{lambda2mu2} we know that $\lambda^{(2)}=\frac{-1+\sqrt{4r-3}}{2}$. Note that \cite[Theorem 2]{DSS} contains a typo in the numerator of the right hand-side of the inequality (the $1$ in the numerator should be a $-1$). The informal description of the result above is that when $\beta$ is large, then the order of the graph will be smaller than the Moore bound. Note however that if $\beta$ is small, then $G_D(\beta)$ may be negative and inequality \eqref{moorebndlarge} may be worse than the classical Moore bound \eqref{moorebnd}. Our results from \cite{CKNV16} may be used to handle some cases when $\beta$ is small (actually when $\lambda_2$ is small). More precisely, Theorem \ref{cknv16thm} gives an upper bound for the order of an $r$-regular graph with small second largest eigenvalue regardless of its diameter actually. 

\begin{table}[h!]
\begin{center}
\begin{tabular}{ |c|c|c|c|c|c| } 
\hline
$(r,D)$ & Known & Defect & Lower & Moore & Upper \\ 
\hline
(8,2) & 57 & 8 & 2.09503 & 2.19258 & 3.40512\\
\hline
(9,2) & 74 & 8 & 2.29956 & 2.37228 & 3.53113\\
\hline
(10,2) & 91 & 10 & 2.46923 & 2.54138 & 3.88473\\
\hline
(4,3) & 41 & 12 & 2.11232 & 2.25342 & 2.88396\\
\hline
(5,3) & 72 & 34 & 2.42905 & 2.62620 & 3.77862\\
\hline
(4,4) & 98 & 63 & 2.53756 & 2.69963 & 3.44307\\
\hline
(5,4) & 212 & 214 & 2.91829 & 3.12941 & 4.41922\\
\hline
(3,5) & 70 & 24 & 2.32340 & 2.39309 & 2.64401\\
\hline 
(4,5) & 364 & 121 & 2.89153 & 2.93996 & 3.42069\\
\hline 
(3,6) & 132 & 58 & 2.45777 & 2.51283 & 2.75001\\
\hline
(4,6) & 740 & 717 & 3.00233 & 3.08314 & 3.73149\\
\hline
\end{tabular}
\end{center}
\caption{Numerical results for small $(r,D)$}
\label{table:1}
\end{table}

We explain this argument and give some numerical examples in Table \ref{table:1} where we listed some values of $(r,D)$ (these values are from \cite[page 4]{MS}) where the maximum orders $n_{r,D}$ of $r$-regular graphs of diameter $D$ are not known. For each such pair, the column labeled {\em Known} gives the largest known order of an $r$-regular graph of diameter $D$. The column {\em Defect} equals the difference between the Moore bound $m_{r,D}$ and the entry in the Known column. The column {\em Moore} contains the value of $\lambda^{(D)}$ rounded below to $5$ decimal points. The column {\em Upper} contains the lower bound for $\tau$ that guarantees that inequality \eqref{moorebndlarge} will give a lower bound than the value from the Known column. For example, for $r=8$ and $D=2$, if $\Ga$ is an $8$-regular graph with diameter $2$ having $\tau<3.40512$, then $|V(\Ga)|<57$. The column {\em Lower} contains an upper bound for $\lambda_2$ that guarantees that the order of such $r$-regular graph would be small. For example, for $r=8$ and $D=2$, our Theorem \ref{cknv16thm} implies that if $\Ga$ is an $8$-regular graph with $\lambda_2<2.0953$, then $|V(\Ga)|<57$. Another way to interpret these results in Table \ref{table:1} is that if one wants to look for a $3$-regular graph of diameter $6$ with more than $132$ vertices, then the second largest eigenvalue of such putative graph has to be between $2.45777$ and $2.75001$.

\section*{\S6 Acknowledgments} This article is based on a talk I gave at the RIMS Conference {\em Research on algebraic combinatorics, related groups and algebras}. I thank the participants of the conference and I am grateful to Hiroshi Nozaki for his help with the computations in Table \ref{table:1} and for his amazing support.

\end{document}